\documentclass[12pt]{article}
\usepackage{amsmath,amssymb,amsfonts,amsthm,url}
\usepackage[margin=1in]{geometry}

\usepackage{graphicx}
\usepackage{subcaption}

\theoremstyle{plain}
\newtheorem{theorem}{Theorem}[section]

\newtheorem{proposition}[theorem]{Proposition}

\theoremstyle{remark}

\newcommand{\NT}{\lfloor n/2 \rfloor}
\newcommand{\Z}{\mathbb{Z}}

\title{On the bigenus of the complete graphs}
\author{Timothy Sun\\San Francisco State University\\Department of Computer Science}
\date{}

\begin{document}

\maketitle

\begin{abstract}
We describe an infinite family of edge-decompositions of complete graphs into two graphs, each of which triangulate the same orientable surface. Previously, such decompositions had only been known for only a few complete graphs. These so-called biembeddings solve a generalization of the Earth-Moon problem for an infinite number of orientable surfaces.
\end{abstract}

\section{Introduction}

Heawood's seminal 1890 paper \cite{Heawood-MapColour} introduced an approach to coloring surface-embedded graphs by upper-bounding their minimum degree. Proving that there are graphs that require as many colors as is used by Heawood's method was accomplished much later by Ringel, Youngs, Jackson, and others \cite{Ringel-MapColor, JacksonRingel-Empire}. Since then, many variations of Heawood's problems have been proposed, most of which have not been solved in general. 

A graph is said to be \emph{biembeddable} in a surface $S$ if it can be expressed as the union of two spanning subgraphs, both embeddable in $S$. Define $\chi_2(S)$ to be the \emph{bichromatic number} of $S$, the maximum chromatic number of any graph biembeddable in $S$. Jackson and Ringel \cite{JacksonRingel-Variations} noted the following variation of Heawood's formula, which they conjectured to be an equality:
\begin{proposition}[Jackson and Ringel \cite{JacksonRingel-Variations}]
Let $S_g$ be the orientable surface of genus $g$, where $g \geq 1$. Then
$$\chi_2(S_g) \leq \left\lfloor \frac{13+\sqrt{73+96g}}{2} \right\rfloor.$$
\label{prop-color}
\end{proposition}
Determining the bichromatic number of the sphere, which is not encompassed by the above formula, is sometimes referred to as the Earth-Moon problem. 

Just like in the Map Color Theorem \cite{Ringel-MapColor}, the most direct approach to showing that this upper bound is tight is through finding biembeddings of complete graphs. Let $\beta(G)$ denote the \emph{bigenus} of $G$, the minimum genus over all orientable surfaces in which $G$ has a biembedding. Cabaniss and Jackson \cite{CabanissJackson-Biembeddings} computed the following lower bound for the complete graphs, which they likewise conjectured to be an equality for sufficiently large $n$:

\begin{proposition}[Cabaniss and Jackson \cite{CabanissJackson-Biembeddings}]
$$\beta(K_n) \geq \left\lceil \frac{n^2-13n+24}{24} \right \rceil.$$
\label{prop-genus}
\end{proposition}

Prior to the present paper, there were only a handful of results relating to the bigenus of the complete graphs. Ringel \cite{Ringel-Farbungsprobleme} showed that $\beta(K_8) = 0$ by decomposing $K_8$ into two planar graphs. While the lower bound suggests that $K_9$ and $K_{10}$ could also have bigenus 0, Battle \emph{et al.} and Tutte \cite{BattleHararyKodama, Tutte-K9} found that $\beta(K_9) > 0$. The statements $\beta(K_{13}) = 1$ and $\beta(K_{14}) = 2$ are special cases of more general results of Ringel \cite{Ringel-Toroidal} and Beineke \cite{Beineke-TwoTorus}, respectively. Cabaniss and Jackson \cite{CabanissJackson-Biembeddings} found triangular biembeddings for $K_{37}$, $K_{61}$, and $K_{85}$ using current graphs, a covering space construction used heavily in the proof of the Map Color Theorem. 

Our main result is a family of current graphs that extends the discoveries of Cabaniss and Jackson to all remaining graphs of the form $K_{24s+13}$, showing for the first time that both conjectures are true at infinitely many values. In Appendix \ref{sec-appendix}, we also provide triangular biembeddings of $K_{16}$, $K_{21}$, and $K_{24}$. 

\section{Graph embeddings}

For background on topological graph theory, see Gross and Tucker \cite{GrossTucker}. For more information about the specific problem and techniques presented here, see Beineke \cite{Beineke-Survey} and Ringel \cite{Ringel-MapColor}. In particular, we assume prior knowledge of current graphs. 

\subsection{Edge bounds from Euler's equation}

A \emph{cellular embedding} of a graph $G$ in the orientable surface of genus $g$ is an injective mapping $\phi: G \to S_g$, where the components of $S_g \setminus \phi(G)$ are open disks. We call these disks \emph{faces}. If the set of faces is denoted by $F(\phi)$, then its size is determined by the \emph{Euler polyhedral equation}
$$|V(G)|-|E(G)|+|F(\phi)| = 2-2g.$$
If $G$ is a simple graph on at least three vertices embedded in $S_g$, then each face has length at least 3. In conjunction with the Euler polyhedral equation, this fact can be used to show that the number of edges in $G$ is at most
$$|E(G)| \leq 3|V(G)|-6+6g,$$
where equality is achieved when all faces are triangular. For biembeddings, both embeddings are bound by the same inequality, so the number of feasible edges doubles:
\begin{proposition}
If $G$ is a simple graph biembedded in the orientable surface $S_g$ of genus $g$, then
$$|E(G)| \leq 6|V(G)|-12+12g,$$
where we have equality if and only if both embeddings are triangular.
\label{prop-edgebound}
\end{proposition}

Propositions \ref{prop-color} and \ref{prop-genus} are straightforward consequences of Proposition \ref{prop-edgebound}, and furthermore, there is a one-way relationship between the two aforementioned conjectures:
\begin{proposition}
If $$\beta(K_n) = \left\lceil \frac{n^2-13n+24}{24} \right \rceil$$ for all $n \geq 13$, then $$\chi_2(S_g) = \left\lfloor \frac{13+\sqrt{73+96g}}{2} \right\rfloor.$$
\end{proposition}
These claims can all be proven by following a nearly identical approach to that in \S4 of Ringel \cite{Ringel-MapColor} used for the original Map Color Theorem. 

When $n \equiv 0, 13, 16, 21 \pmod{24}$, the expression inside the ceiling function in the bigenus lower bound is an integer, so any biembedding of $K_n$ achieving this lower bound must consist of two triangular embeddings. For brevity, we also use the word \emph{triangular} to describe such a biembedding. All of our constructions will be triangular biembeddings. 

\subsection{Rotation systems}

To describe an embedding combinatorially, we associate, with each edge $e \in E$, two arcs $e^+$ and $e^-$ with the same endpoints, each representing the two different directions in which $e$ can be traversed. The set of such arcs is denoted $E^+$. A \emph{rotation} of a vertex is a cyclic permutation of the arcs leaving that vertex, and a \emph{rotation system} of a graph is an assignment of a rotation to each vertex. The Heffter-Edmonds principle (see \S3.2 of Gross and Tucker \cite{GrossTucker}) states that rotation systems are in one-to-one correspondence with cellular embeddings in orientable surfaces, up to orientation-preserving equivalence of embeddings. The cyclic orientation of the edges incident with each vertex induces a rotation system, and each rotation system can be converted into an embedding by face-tracing.

\subsection{Circulant and current graphs}

Let $\Z_n$ denote the integers modulo $n$, and let $X$ be any set of nonzero elements of $\Z_n$. The \emph{circulant graph} $C(n, X)$ is the simple graph with vertex set $\Z_n$, where the neighbors of each vertex $i$ are the vertices $i \pm x$ for each $x \in X$. The complete graph $K_n$ is the circulant graph $C(n, \{1, 2, \dotsc, \NT\})$, and we take the convention that $X$ is always a subset of $\{1, 2, \dotsc, \NT\}$. One way of decomposing the edges of the complete graph would be to find a partition of $\{1, 2, \dotsc, \lfloor n/2 \rfloor\}$ into two equal-sized subsets $X_1$ and $X_2$ and to form the circulant graphs $C(n, X_1)$ and $C(n, X_2)$. Since circulant graphs have cyclic symmetry, triangular embeddings of such graphs can potentially be constructed using the theory of current graphs.

A \emph{current graph} is an arc-labeled, embedded graph, where, for our purposes, the graph is embedded in an orientable surface, and the arc-labeling $\alpha: E^+ \to \Z_n \setminus \{0\}$ satisfies $\alpha(e^+) = -\alpha(e^-)$ for each edge $e$. We call $\Z_n$ the \emph{current group} and the arc labels \emph{currents}. Our current graphs also satisfy some standard properties:

\begin{itemize}
\item The embedding has one face. 
\item Each vertex has degree 3.
\item The currents entering each vertex sum to 0. 
\item Each nonzero element in $\Z_n$ appears at most once as a current.
\end{itemize}

If all of these properties are satisfied and the currents are $\{\pm c_1, \pm c_2, \dotsc, \pm c_i\}$, then the derived embedding is an orientable triangular embedding of the circulant graph $$C(n, \{c_1, c_2, \dotsc, c_i\})$$ (see \S2.3 of Ringel \cite{Ringel-MapColor} or \S4 of Gross and Tucker \cite{GrossTucker}). We draw our current graphs so that each edge is represented by the arc with the ``smaller'' current, conforming to our convention for circulant graphs. 

\begin{figure}[!t]
\centering
\includegraphics[scale=1.1]{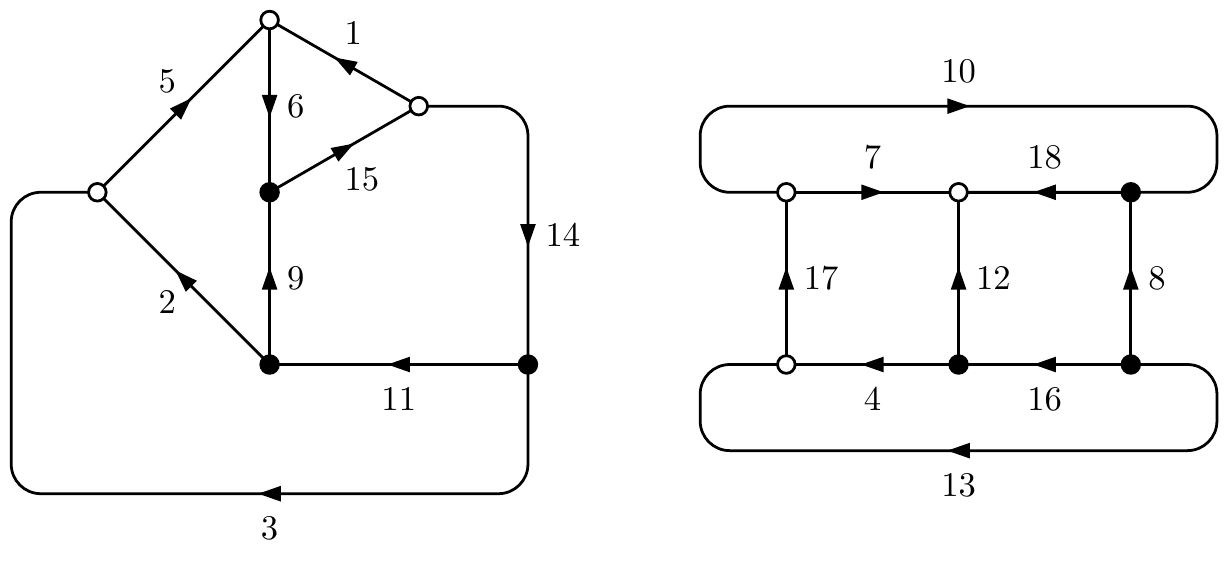}
\caption{A pair of current graphs with current group $\mathbb{Z}_{37}$, where solid and hollow vertices represent clockwise and counterclockwise rotations, respectively.}
\label{fig-bi1}
\end{figure}

In our example in Figure \ref{fig-bi1}, the group elements $1, 2, \dotsc, 18$ are partitioned evenly between the two current graphs, so the two derived embeddings form a triangular biembedding of $K_{37}$. In general, pairs of evenly-sized current graphs satisfying the aforementioned properties can be used to find triangular biembeddings of $K_{24s+13}$, for $s \geq 0$. Such pairs were previously only known for $s \leq 3$ \cite{Ringel-Toroidal, CabanissJackson-Biembeddings}. 

\section{The main construction}

\begin{figure}[!t]
\centering
\includegraphics[scale=1.1]{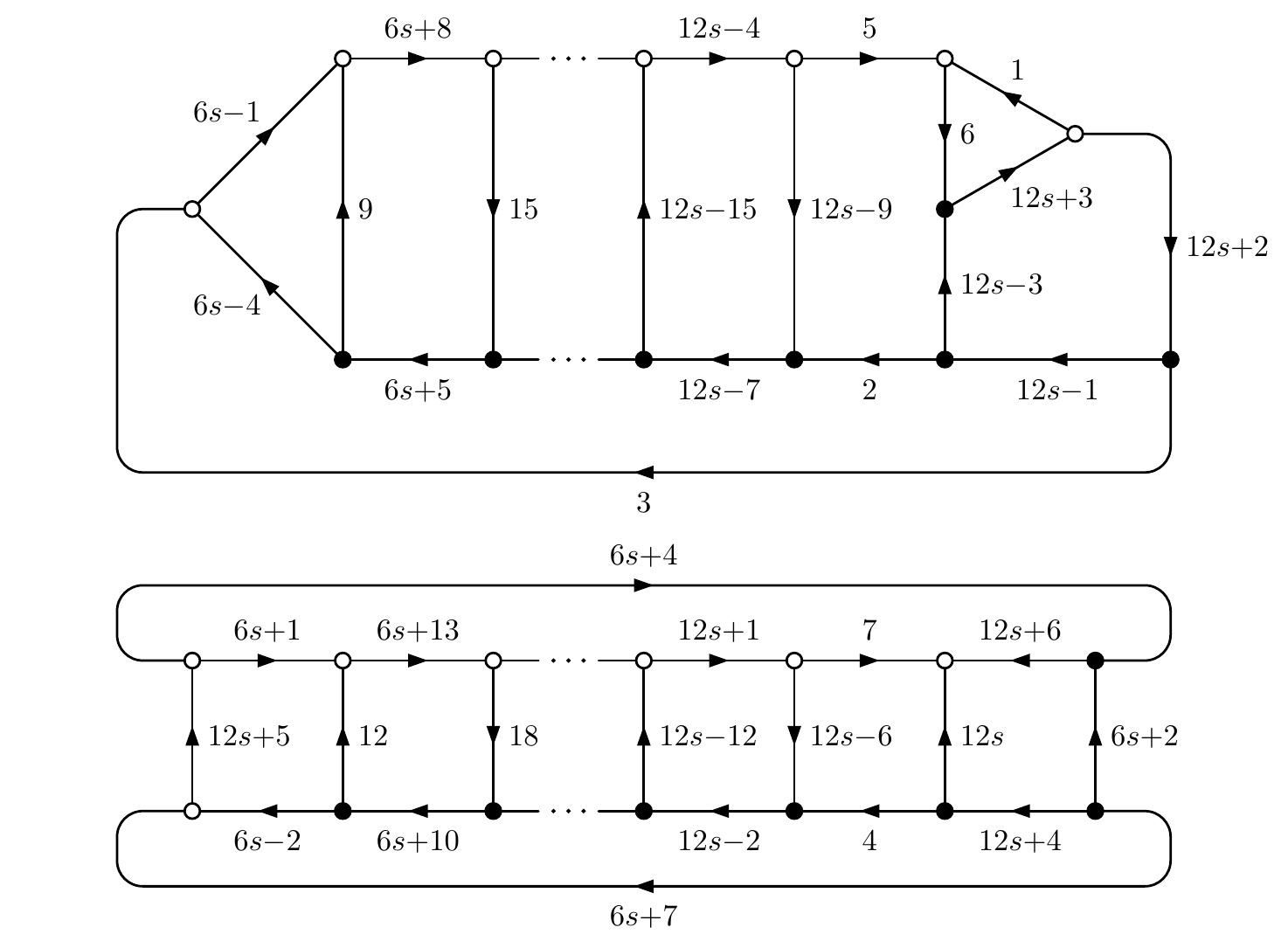}
\caption{The general construction with current group $\mathbb{Z}_{24s+13}$, $s \geq 1$.}
\label{fig-bi2}
\end{figure}

The family of current graphs in Figure \ref{fig-bi2} generates triangular biembeddings of the complete graphs $K_{24s+13}$, for all $s \geq 1$. The current graphs previously shown in Figure \ref{fig-bi1} constitute the $s=1$ case, the smallest pair in the family. Both current graphs contain a ``ladder'' whose ``rungs'' have currents that alternate in direction and form an arithmetic sequence with step size 6. For each $s$, the labels $1, \dotsc, 12s{+}6$ are partitioned between the two graphs: besides the exceptional labels $1$, $6$, $6s{+}2$, and $12s{+}5$, the first current graph in each pair has all the labels $i$ such that $i \equiv 2,3,5 \pmod{6}$, and the second current graph has all the labels $j$ such that $j \equiv 0, 1, 4 \pmod{6}$. We make one minor but necessary check to ensure that the derived graphs are connected. Otherwise, the embedding could consist of more than one triangulated surface. Our current graphs always produce connected graphs because they have a generator of the current group (e.g. $1$ or $4$) as a current and hence have Hamiltonian cycles.

\bibliographystyle{alpha}
\bibliography{biblio}

\newpage
\appendix

\section{Self-complementary biembeddings}\label{sec-appendix}

Given a simple graph $G$, its \emph{edge-complement} $\overline{G}$ is the simple graph on the same vertex set $V(G)$, where for every pair of vertices $u,v \in V(G)$, $(u,v)$ is an edge of $\overline{G}$ if and only if $(u,v)$ is not an edge of $G$. A \emph{self-complementary} graph $G$ is one which is isomorphic to its own edge-complement $\overline{G}$, and an isomorphism from $G$ to $\overline{G}$ is known as an \emph{antimorphism}. We note that since the edge-complement of a disconnected graph is connected, all self-complementary graphs are connected. 

A triangular embedding of a self-complementary graph on $n$ vertices automatically leads to a triangular biembedding of the complete graph $K_n$ by simply reusing the same embedding on its edge-complement. The rotation systems given in Tables \ref{tab-sc16}, \ref{tab-sc21}, and \ref{tab-sc24} describe triangular embeddings of self-complementary graphs on $16$, $21$, and $24$ vertices found through computer search.

\begin{table}[!ht]
\centering
{\footnotesize
$$\begin{array}{rrrrrrrrrrrrrrrrrrrrrrrrrr}
0. & 1 & 9 & 5 & 3 \\
1. & 0 & 3 & 11 & 5 & 12 & 15 & 14 & 7 & 13 & 8 & 9 \\
2. & 3 & 5 & 11 & 7 \\
3. & 0 & 5 & 2 & 7 & 14 & 15 & 10 & 13 & 9 & 11 & 1 \\
4. & 5 & 7 & 9 & 13 \\
5. & 0 & 9 & 15 & 7 & 4 & 13 & 12 & 1 & 11 & 2 & 3 \\
6. & 7 & 11 & 15 & 9 \\
7. & 1 & 14 & 3 & 2 & 11 & 6 & 9 & 4 & 5 & 15 & 13 \\
8. & 1 & 13 & 11 & 9 \\
9. & 0 & 1 & 8 & 11 & 3 & 13 & 4 & 7 & 6 & 15 & 5 \\
10. & 3 & 15 & 11 & 13 \\
11. & 1 & 3 & 9 & 8 & 13 & 10 & 15 & 6 & 7 & 2 & 5 \\
12. & 1 & 5 & 13 & 15 \\
13. & 1 & 7 & 15 & 12 & 5 & 4 & 9 & 3 & 10 & 11 & 8 \\
14. & 1 & 15 & 3 & 7 \\
15. & 1 & 12 & 13 & 7 & 5 & 9 & 6 & 11 & 10 & 3 & 14
\end{array}$$
}
\caption{A triangular embedding of a self-complementary graph on 16 vertices.}
\label{tab-sc16}
\end{table}

\begin{table}[!ht]
\centering
{\footnotesize
$$\begin{array}{rrrrrrrrrrrrrrrrrrrrrrrrrr}
0. & 1 & 7 & 9 & 5 & 3 \\
1. & 0 & 3 & 11 & 20 & 9 & 17 & 16 & 5 & 13 & 12 & 15 & 14 & 19 & 18 & 7 \\
2. & 3 & 7 & 11 & 5 & 9 \\
3. & 0 & 5 & 18 & 19 & 16 & 17 & 14 & 15 & 7 & 2 & 9 & 20 & 13 & 11 & 1 \\
4. & 5 & 11 & 9 & 7 & 13 \\
5. & 0 & 9 & 2 & 11 & 4 & 13 & 1 & 16 & 19 & 20 & 15 & 17 & 7 & 18 & 3 \\
6. & 7 & 15 & 9 & 11 & 13 \\
7. & 0 & 1 & 18 & 5 & 17 & 20 & 19 & 11 & 2 & 3 & 15 & 6 & 13 & 4 & 9 \\
8. & 9 & 13 & 15 & 11 & 17 \\
9. & 0 & 7 & 4 & 11 & 6 & 15 & 19 & 13 & 8 & 17 & 1 & 20 & 3 & 2 & 5 \\
10. & 11 & 19 & 15 & 13 & 17 \\
11. & 1 & 3 & 13 & 6 & 9 & 4 & 5 & 2 & 7 & 19 & 10 & 17 & 8 & 15 & 20 \\
12. & 1 & 13 & 19 & 17 & 15 \\
13. & 1 & 5 & 4 & 7 & 6 & 11 & 3 & 20 & 17 & 10 & 15 & 8 & 9 & 19 & 12 \\
14. & 1 & 15 & 3 & 17 & 19 \\
15. & 1 & 12 & 17 & 5 & 20 & 11 & 8 & 13 & 10 & 19 & 9 & 6 & 7 & 3 & 14 \\
16. & 1 & 17 & 3 & 19 & 5 \\
17. & 1 & 9 & 8 & 11 & 10 & 13 & 20 & 7 & 5 & 15 & 12 & 19 & 14 & 3 & 16 \\
18. & 1 & 19 & 3 & 5 & 7 \\
19. & 1 & 14 & 17 & 12 & 13 & 9 & 15 & 10 & 11 & 7 & 20 & 5 & 16 & 3 & 18 \\
20. & 1 & 11 & 15 & 5 & 19 & 7 & 17 & 13 & 3 & 9
\end{array}$$
}
\caption{A triangular embedding of a self-complementary graph on 21 vertices.}
\label{tab-sc21}
\end{table}

\begin{table}[!ht]
\centering
{\footnotesize
$$\begin{array}{rrrrrrrrrrrrrrrrrrrrrrrrrr}
0. & 1 & 13 & 16 & 6 & 21 & 2 & 17 & 20 & 15 & 18 & 8 & 14 & 22 & 19 & 4 & 10 \\
1. & 0 & 10 & 8 & 6 & 4 & 12 & 13 \\
2. & 0 & 21 & 4 & 19 & 18 & 15 & 3 & 10 & 12 & 20 & 16 & 8 & 23 & 6 & 22 & 17 \\
3. & 2 & 15 & 14 & 6 & 12 & 8 & 10 \\
4. & 0 & 19 & 2 & 21 & 20 & 14 & 8 & 5 & 17 & 22 & 18 & 12 & 1 & 6 & 23 & 10 \\
5. & 4 & 8 & 12 & 14 & 10 & 16 & 17 \\
6. & 0 & 16 & 10 & 20 & 12 & 3 & 14 & 7 & 19 & 22 & 2 & 23 & 4 & 1 & 8 & 21 \\
7. & 6 & 14 & 16 & 12 & 10 & 18 & 19 \\
8. & 0 & 18 & 22 & 23 & 2 & 16 & 9 & 21 & 6 & 1 & 10 & 3 & 12 & 5 & 4 & 14 \\
9. & 8 & 16 & 14 & 12 & 18 & 20 & 21 \\
10. & 0 & 4 & 23 & 11 & 20 & 6 & 16 & 5 & 14 & 18 & 7 & 12 & 2 & 3 & 8 & 1 \\
11. & 10 & 23 & 22 & 16 & 18 & 14 & 20 \\
12. & 1 & 4 & 18 & 9 & 14 & 5 & 8 & 3 & 6 & 20 & 2 & 10 & 7 & 16 & 22 & 13 \\
13. & 0 & 1 & 12 & 22 & 20 & 18 & 16 \\
14. & 0 & 8 & 4 & 20 & 11 & 18 & 10 & 5 & 12 & 9 & 16 & 7 & 6 & 3 & 15 & 22 \\
15. & 0 & 20 & 22 & 14 & 3 & 2 & 18 \\
16. & 0 & 13 & 18 & 11 & 22 & 12 & 7 & 14 & 9 & 8 & 2 & 20 & 17 & 5 & 10 & 6 \\
17. & 0 & 2 & 22 & 4 & 5 & 16 & 20 \\
18. & 0 & 15 & 2 & 19 & 7 & 10 & 14 & 11 & 16 & 13 & 20 & 9 & 12 & 4 & 22 & 8 \\
19. & 0 & 22 & 6 & 7 & 18 & 2 & 4 \\
20. & 0 & 17 & 16 & 2 & 12 & 6 & 10 & 11 & 14 & 4 & 21 & 9 & 18 & 13 & 22 & 15 \\
21. & 0 & 6 & 8 & 9 & 20 & 4 & 2 \\
22. & 0 & 14 & 15 & 20 & 13 & 12 & 16 & 11 & 23 & 8 & 18 & 4 & 17 & 2 & 6 & 19 \\
23. & 2 & 8 & 22 & 11 & 10 & 4 & 6
\end{array}$$
}
\caption{A triangular embedding of a self-complementary graph on 24 vertices.}
\label{tab-sc24}
\end{table}

These self-complementary graphs were constructed starting from antimorphisms with cycle decompositions
$$(0 \,\,\,\, 1 \,\,\,\, 2 \,\,\,\, \dotsc \,\,\,\, n-1)$$
for $n = 16\text{~and~}24$, and 
$$(0 \,\,\,\, 1 \,\,\,\, 2 \,\,\,\, \dotsc \,\,\,\, n-2)(n-1)$$
for $n = 21$. When the antimorphism takes on one of these two forms, the self-complementary graph is uniquely determined by the edges leaving, for example, vertex 0. The number of self-complementary graphs grows quite rapidly, hence it seems very likely that a triangular embedding for some self-complementary graph on $n$ vertices exists for all $n \equiv 0, 13, 16, 21 \pmod{24}$. We note that Ringel's biembeddings of $K_8$ and $K_{13}$ \cite{Ringel-Farbungsprobleme, Ringel-Toroidal} were derived from self-complementary graphs, as well.

\end{document}